\documentstyle[12pt]{article}
\begin{document}
\title{A generalization of Kostant theorem to integral cohomology}
\author {Qibing Zheng\\School of Mathematical Science and LPMC, Nankai University\\
Tianjin 300071, China\\
zhengqb@nankai.edu.cn\footnote{Project Supported by Natural Science
Foundation of China, grant No. 11071125}} \maketitle
\input amssym.def
\newsymbol\leqslant 1336
\newsymbol\geqslant 133E
\baselineskip=20pt
\def\w{\widetilde}
\def\o{\overline}

\begin{abstract} In this paper, we find weight decomposition and rank of a weight
in the integral (co)homology of the positive system of a semi-simple Lie algebra
over $\Bbb C$ and prove that the (co)homology
of the weight subcomplex over a field of characteristic $p$ is $0$ if the rank of the weight
is not divisible by $p$.
This generalizes Kostant theorem to the integral cohomology of the positive system.
\end{abstract}
\section*{\hspace*{40mm}Table of Contents}

\hspace*{19.3mm}section 1 Introduction

\hspace*{13mm}section 2 Basis graph, weight and Weyl group action

\hspace*{13mm}section 3 Weight subcomplex

\hspace*{13mm}section 4 Index and neighbors

\hspace*{13mm}section 5 Basic properties of basis graph

\hspace*{13mm}section 6 Diamond (co)chain complex

\hspace*{13mm}section 7 Combinatorial weight

\newtheorem{Definition}{Definition}[section]
\newtheorem{Theorem}{Theorem}[section]
\newtheorem{Lemma}{Lemma}[section]
\newtheorem{Example}{Example}[section]
\font\hua=eusm10 scaled\magstephalf

\section{Introduction}\vspace{3mm}

{\hspace*{5.5mm} The cohomology of Lie algebras is a very important topic in many fields of Mathematics.
Kostant in \cite{k} compute ${\rm Ext}_{\Bbb C(\frak R^+)}^*(\Bbb C,V)$, where $V$ is a representation of a semi-simple Lie algebra $\frak G$ over $\Bbb C$ and $\Bbb C(R^+)$ is the Lie algebra over $\Bbb C$ generated
by the positive root system $R^+$ of $\frak G$.
As a special case, $H^*(R^+;\Bbb C)={\rm Ext}_{\Bbb C(R^+)}^*(\Bbb C,\Bbb C)$ has a basis in 1-1 correspondence
with the Weyl group of $\frak G$.
But very little is known about the torsion part of $H^*(R^+;\Bbb Z)={\rm Ext}^*_{\Bbb Z(R^+)}(\Bbb Z,\Bbb Z)$.
This problem is related to homotopy theory.
Let $A^+_n$ be the positive root system of the simple Lie algebra over $\Bbb C$ with Dynkin graph $A_n$
and $A^+_{\infty}=\cup_{n}\,A^+_n$.
Then $H^*(A^+_{\infty};\Bbb Z_p)=\cup_n\,H^*(A^+_{n};\Bbb Z_p)$
($\Bbb Z_p$ is the field of integers modular a prime $p$) is
a direct sum summand of May spectral sequence \cite{m} converging to the cohomology of the
Steenrod algebra (see \cite{d}). The cohomology of Steenrod algebra is the $E_2$-term of the Adams
spectral sequence converging to the stable homotopy groups of spheres (see \cite{a} and \cite{r}).
In \cite{h}, we find weight and rank in the integral cohomology  $H^*(A^+_n;\Bbb Z)$
and prove a lot of properties of the cohomology that includes the work of Dwyer in \cite{v}.
In this paper, we generalize the result in \cite{h} to all positive root system of a semi-simple
Lie algebra over $\Bbb C$.

In section 2, we first define the basis graph $\Gamma(R^+)$ of a positive root system
$R^+$ of a semi-simple Lie algebra over $\Bbb C$ in Definition 2.2. Then we define Weyl group
action $\circ$ on $\Gamma(R^+)$ in Theorem 2.3.
In Theorem 2.4, we define weight on $\Gamma(R^+)$ and find weight decomposition
$\Gamma(R^+)=\sqcup_{\alpha\in\Omega(R^+)}\,\Gamma(\alpha)$.

In section 3, we prove in Theorem 3.2 the weight complex decomposition
$$(\Lambda(R^+),d)=\oplus_{\alpha\in\Omega(R^+)}(\Lambda(\alpha),d),\quad
(\Lambda(R^+),\delta)=\oplus_{\alpha\in\Omega(R^+)}(\Lambda(\alpha),\delta),$$
that is induced by the weight decomposition of basis graph and so have the (co)homology weight decomposition
$$H_*(R^+)=\oplus_{\alpha\in\Omega(R^+)}H_*^\alpha(R^+),\quad
H^*(R^+)=\oplus_{\alpha\in\Omega(R^+)}H^*_\alpha(R^+).$$
In Theorem 3.3, we prove the duality theorem on weight subcomplexes.

In section 4, we define the index of a weight with respect of a positive root in Definition 4.2
and find the relation between index and neighbors in Theorem 4.4 that is very important to
the proof of properties of basis graph in section 5.

In section 5, we prove the basic properties of basis graph.
In Definition 5.1 and Theorem 5.2, we prove that basis graph $\Gamma(R^+)$ is a diamond graph.
In Theorem 5.3, we prove that the weight decomposition of $\Gamma(R^+)$ coincides with its connected component
decomposition.
In Definition 5.4 and Theorem 5.5, we prove that $\Gamma(R^+)$ is an admissible diamond graph.

In section 6, we prove the property of $H^*(R^+;\Bbb F_p)$ with $\Bbb F_p$
a field of characteristic $p$.
In Definition 6.1 and Theorem 6.2, we prove that $(\Lambda(R^+),\delta)$ is a diamond cochain complex.
In Definition 6.3 and Theorem 6.4, we prove that for a diamond cochain complex $C$ with basis graph $\Gamma$,
$H^*(C;\Bbb F_p)=0$ if the basis graph $\Gamma$ is connected with rank
not divisible by $p$. As a corollary of Theorem 6.4, we compute in Theorem 6.7
the free part of the integral cohomology $H^*(R^+;\Bbb Z)$ and prove that $H^*_{\alpha}(R^+;\Bbb F_p)=0$
if $r(\alpha)$ is not divisible by the prime $p$, i.e.,
$$H_*^{free}(R^+)=\oplus_{w\in W}\,H_*^{w{*}\varrho}(R^+),\,\,
H^*_{free}(R^+)=\oplus_{w\in W}\,H^*_{w{*}\varrho}(R^+).$$
$$H_*(R^+;\Bbb F_p)=\oplus_{p|r(\alpha)}\,H_*^\alpha(R^+;\Bbb F_p),\,\,
H^*(R^+;\Bbb F_p)=\oplus_{p|r(\alpha)}\,H^*_\alpha(R^+;\Bbb F_p).$$
This generalizes Kostant theorem to the integral cohomology.

In section 7, we compare the combinatorial weight defined in \cite{h} and weight in this paper
for the positive root system of the simple Lie algebra with Dynkin graph $A_n$.

\section{Basis graph, weight and Weyl group action}\vspace{3mm}

\hspace*{5.5mm}{\bf Definition 2.1} In this paper,
$R^+$ is the positive root system of a semi-simple Lie algebra
over the complex field $\Bbb C$.
$W$ is the Weyl group acting on the real vector space $V$ spanned by $R^+$.
A W-invariant metric is given on $V$.

We denote the action of $W$ on $V$ by ${*}$, i.e., for $w\in W$ and $v\in V$,
the action of $w$ on $v$ is denoted by $w{*}v$.
For $e\in R^+$, $r_e\in W$ is the reflection on $V$ keeping the hyperplane orthogonal to $e$ fixed.
\vspace{3mm}

{\bf Conventions}
In this paper, we use the same symbol to denote both a graph (finite 1-dimensional simplicial complex)
and its vertex set.
So for a graph $\Gamma$ and two vertices $a,b\in\Gamma$, $\{a,b\}\in\Gamma$ implies $\{a,b\}$ is an edge of $\Gamma$
and $\{a,b\}\not\in\Gamma$ implies $\{a,b\}$ is not an edge of $\Gamma$.\vspace{3mm}

{\bf Definition 2.2} The basis graph $\Gamma(R^+)$ is defined as follows.
The vertex set $\Gamma(R^+)$ is the set of all subsets of $R^+$ (including empty set $\emptyset$).
For two vertices $\sigma,\tau\in\Gamma(R^+)$,
$\{\sigma,\tau\}\in\Gamma(R^+)$ if there are $e,e_1,e_2\in R^+$ such that $e=e_1{+}e_2$ and\vspace{1mm}\\
\hspace*{9.7mm}$\sigma{\setminus}\tau=\{e\}\,\,{\rm and}\,\,\tau{\setminus}\sigma=\{e_1,e_2\},\,\,{\rm or}\,\, \tau{\setminus}\sigma=\{e\}\,\,{\rm and} \,\,\sigma{\setminus}\tau=\{e_1,e_2\}$.
\vspace{3mm}

{\bf Theorem 2.3} {\it
There is a Weyl group action on the vertex set of the basis graph (denoted by ${\circ}$)
$${\circ}\,\colon W{\times}\Gamma(R^+)\to \Gamma(R^+)$$
defined as follows.
For $w\in W$ and a subset $S$ of $V$, let $w{*}S=\{w{*}v\,|\,v\in S\}$.
Then for $\sigma\in\Gamma(R^+)$ ($\sigma^c=R^+{\setminus}\sigma$),\vspace{2mm}\\
\hspace*{45mm}$w{\circ}\sigma=w{*}(\sigma{\cup}{-}\sigma^c)\cap R^+$.
}\vspace{3mm}

{\it Proof}\, Suppose $w,w'\in W$. Then for any $\sigma\subset R^+$,
$$R^+=\Big((w{*}\sigma){\cup}(w{*}{-}\sigma){\cup}(w{*}\sigma^c){\cup}(w{*}{-}\sigma^c)\Big){\cap}R^+.$$
So ($R^-=-R^+=\{-e\,|\,e\in R^+\}$)
\begin{eqnarray*}&&w'{\circ}(w{\circ}\sigma)\\
&=&w'{\circ}\Big(w{*}(\sigma{\cup}{-}\sigma^c)\cap R^+\Big)\\
&=&w'{*}\Big((w{*}(\sigma{\cup}{-}\sigma^c)\cap R^+)\cup
{-}(w{*}(\sigma{\cup}{-}\sigma^c)\cap R^+)^c\Big)\cap R^+\\
&=&w'{*}\Big((w{*}(\sigma{\cup}{-}\sigma^c)\cap R^+)\cup
(w{*}(\sigma{\cup}{-}\sigma^c)\cap R^-)\Big)\cap R^+\\
&=&w'{*}(w{*}(\sigma{\cup}{-}\sigma^c))\cap R^+\\
&=&(w'w){*}(\sigma{\cup}{-}\sigma^c)\cap R^+\\
&=&(w'w){\circ}(\sigma)
\end{eqnarray*}
and $1{\circ}\sigma=(\sigma{\cap}R^+){\cup}({-}\sigma^c{\cap}R^+)=\sigma$.
The group action is well-defined.
\hfill$\Box$\vspace{3mm}

{\bf Remark} For a simple reflection $r_e$, the ${*}$ action of $r_e$ on $R^+{\setminus}\{e\}$ is a permutation
and the $\circ$ action of $r_e$ on $\Gamma(R^+)$ is as follows
$$\hspace{31mm}r_e{\circ}\sigma=\left\{
\begin{array}{ccc}
r_e{*}(\sigma{\setminus}\{e\})&&{\rm if}\,\,e\in\sigma,\\
(r_e{*}\sigma)\cup\{e\}&&{\rm if}\,\,e\not\in\sigma.
\end{array}
\right.\hspace{31mm}(2.1)$$
But the group action $\circ$ is not free, i.e., the isotropy group of a vertex may be non-trivial.
\vspace{3mm}

{\bf Theorem 2.4}\, {\it There is a weight set $\Omega(R^+)\subset V$ containing the half sum of positive roots
$\varrho=\frac 12\Sigma_{e\in R^+}e$ such that
$$\Gamma(R^+)=\sqcup_{\alpha\in \Omega(R^+)}\,\Gamma(\alpha),$$
where $\Gamma(\alpha)$ is a full subgraph (called weight subgraph) of $\Gamma(R^+)$
and $\sqcup$ means the disjoint union of graphs.

The $\circ$ group action in Theorem 2.3 keeps the weight and graph structure, i.e.,
for any $w\in W$,  the 1-1 correspondence $w{\circ}\,\colon \Gamma(R^+)\to \Gamma(R^+)$ is a 
graph isomorphism such that the restriction on $\Gamma(\alpha)$ is a graph isomorphism 
from $\Gamma(\alpha)$ to $\Gamma(w{*}\alpha)$ for all $\alpha\in\Omega(R^+)$.
}\vspace{3mm}

{\it Proof}\, Define weight map $\omega\colon\Gamma(R^+)\to V$ as follows.
For $\sigma\in\Gamma(R^+)$, $\omega(\sigma)=\varrho-\Sigma_{e\in \sigma}\,e$.
Define $\Omega(R^+)={\rm im}\,\omega$. Then $\varrho=\omega(\emptyset)\in\Omega(R^+)$.
For a weight $\alpha\in\Omega(R^+)$, define $\Gamma(\alpha)=\omega^{-1}(\alpha)$.
It is obvious that if $\{\sigma,\tau\}\in\Gamma(R^+)$, then $\omega(\sigma)=\omega(\tau)$.
So $\Gamma(\alpha)$ is a full subgraph of $\Gamma(R^+)$ for all $\alpha\in\Omega(R^+)$
and $\Gamma(R^+)=\sqcup_{\alpha\in\Omega(R^+)}\,\Gamma(\alpha)$.

For any $w\in W$ and $\sigma\in\Gamma(R^+)$,\\
\hspace*{34mm}$w{*}(\omega(\sigma))$\vspace{1mm}\\
\hspace*{30mm}$=w{*}(\varrho-\Sigma_{e\in \sigma}\,e)$\vspace{1mm}\\
\hspace*{30mm}$=w{*}\varrho-\Sigma_{e\in \sigma}\,w{*}e$\vspace{1mm}\\
\hspace*{30mm}$=\frac 12(\Sigma_{e\in R^+\!,\,w{*}e\in R^+}\,w{*}e+
\Sigma_{e\in R^+\!,\,w{*}e\in R^-}\,w{*}e)-\Sigma_{e\in \sigma}\,w{*}e$\vspace{1mm}\\
\hspace*{30mm}$=\varrho+\Sigma_{e\in R^+\!,\,w{*}e\in R^-}\,w{*}e-\Sigma_{e\in \sigma}\,w{*}e$\vspace{1mm}\\
\hspace*{30mm}$=\varrho+\Sigma_{e\in \sigma^c,\,w{*}e\in R^-}\,w{*}e-\Sigma_{e\in \sigma,\,w{*}e\in R^+}\,w{*}e$\vspace{1mm}\\
\hspace*{30mm}$=\varrho-\Sigma_{e\in w{\circ}\sigma}\,e$\vspace{1mm}\\
\hspace*{30mm}$=\omega(w{\circ}\sigma).$\vspace{1mm}

So $\sigma\to w{\circ}\sigma$ is a 1-1 correspondence from $\Gamma(\alpha)$ to $\Gamma(w{*}\alpha)$.

Now we prove the map $\sigma\to w{\circ}\sigma$ is a graph isomorphism.
Since $W$ is generated by simple reflections, we need only prove the special case $w$ is a simple
reflection $r_e$.

Suppose $\{\sigma,\tau\}\in\Gamma(\alpha)$, $\sigma{\setminus}\tau=\{f\}$, $\tau{\setminus}\sigma=\{f',f''\}$,
$f=f'{+}f''$ and $e$ is a simple root.
If $e\not\in\{f',f''\}$, then by (2.1),
$$r_e{\circ}\sigma{\setminus}r_e{\circ}\tau=\{r_e{*}f\},\,
r_e{\circ}\tau{\setminus}r_e{\circ}\sigma=\{r_e{*}f',r_e{*}f''\},\,
r_e{*}f=r_e{*}f'{+}r_e{*}f''.$$
So $\{r_e{\circ}\sigma,r_e{\circ}\tau\}\in\Gamma(r_e{*}\alpha)$.
If $e\in\{f',f''\}$, we may suppose $e=f'$, then
$$r_e{\circ}\sigma{\setminus}r_e{\circ}\tau=\{r_e{*}f,e\},\,
r_e{\circ}\tau{\setminus}r_e{\circ}\sigma=\{r_e{*}f''\},\,
r_e{*}f''=r_e{*}f{+}e.$$
So $\{r_e{\circ}\sigma,r_e{\circ}\tau\}\in\Gamma(r_e{*}\alpha)$.
The map $\sigma\to r_e{\circ}\sigma$ is a graph isomorphism.
\hfill$\Box$\vspace{3mm}

{\bf Theorem 2.5} {\it For $\alpha\in\Omega(R^+)$, the following conditions are equivalent.

\begin{enumerate}

\item $\Gamma(\alpha)$ has only one vertex.

\item $\alpha=w{*}\varrho$ for some $w\in W$, equivalently, $\alpha$ is in the orbit of $\varrho$
with respect to the ${*}$ group action of $W$ on $\Omega(R^+)$.

\item The isotropy group of $\alpha$ with respect to the $*$ action of $W$ on $\Omega(R^+)$ is trivial.

\item $r_e{*}\alpha\neq\alpha$ for all reflections $r_e\in W$ with $e\in R^+$.
\end{enumerate}
}\vspace{3mm}

{\it Proof}\, (1)$\Rightarrow$(2). Suppose $\Gamma(\alpha)$ has only one vertex $\sigma$.
If $\sigma=\emptyset$, then $\alpha=\varrho$ and the conclusion holds. Suppose $\sigma\neq\emptyset$.
We first prove a lemma that there must be a simple root in $\sigma$.
Suppose there is no simple root in $\sigma$.
Define partial order on $\sigma$ as follows.
For $e_1,e_2\in \sigma$, $e_1<e_2$ if there is a simple root $s\in R^+$ such that $e_1{+}s=e_2$.
Then there is a smallest root $e\in \sigma$.
Since $e$ is not simple, there is simple root $s_1$
and $e'\not\in \sigma$ such that $e=s_1{+}e'$. So $\sigma'=(\sigma{\setminus}\{e\})\cup\{s_1,e'\}\in\Gamma(\alpha)$.
A contradiction!
So the lemma holds.

Let $e_1$ be the simple root in $\sigma$ and $\alpha_1=r_{e_1}{*}\alpha$.
Since the graphs $\Gamma(\alpha)$ and $\Gamma(\alpha_1)$ are isomorphic,
$\Gamma(\alpha_1)$ has only one vertex $\sigma_1=r_{e_1}{*}(\sigma{\setminus}\{e_1\})$ by (2.1).
By the lemma, there is a simple root $e_{2}\in S_1$.
Let $\alpha_2=r_{e_2}{*}\alpha_1$.
Then $\Gamma(\alpha_2)$ has only one vertex $\sigma_2=r_{e_2}{*}(\sigma_1{\setminus}\{e_2\})$.
Repeat this process and we get $\sigma_1,\sigma_2,\cdots$ such that $|\sigma_{i+1}|=|\sigma_i|{-}1$ ($|\cdot|$ the cardinality).
So there is an $n$ such that $\sigma_n=\emptyset$ and $\alpha_n=(r_{e_n}\cdots r_{e_1}){*}\alpha=\omega(\emptyset)=\varrho$.
$\alpha=(r_{e_1}\cdots r_{e_n}){*}\varrho$.

(2)$\Rightarrow$(3). Since $W$ acts freely on the orbit of $\varrho$,
the isotropy group of $\alpha$ is trivial.

(3)$\Rightarrow$(4). The isotropy group of $\alpha$ is not trivial
if and only if $\alpha$ is in the wall of a Weyl chamber and the isotropy group of $\alpha$ is
generated by the reflections keeping the wall containing $\alpha$ fixed, i.e. $r_e{*}\alpha=\alpha$ for some reflection $r_e\in W$.
So (3) is equivalent to (4).

(4)$\Rightarrow$(1). If $r_e{*}\alpha=\alpha$, let
$$\Gamma_e^1(\alpha)=\{\sigma\in\Gamma(\alpha)\,|\,e\in \sigma\},\,
\Gamma_e^0(\alpha)=\{\sigma\in\Gamma(\alpha)\,|\,e\not\in \sigma\},$$
then $\sigma\to r_e{\circ}\sigma$ is an 1-1 correspondence between $\Gamma_e^1(\alpha)$ and $\Gamma_e^0(\alpha)$.
So $\Gamma(\alpha)$ has more than one vertices.
\vspace{3mm}

\section{Weight subcomplex}\vspace{3mm}

\hspace*{5.5mm}{\bf Definition 3.1}
$\Lambda(R^+)$ is the exterior algebra over $\Bbb Z$ generated by $R^+$.
We have to distinguish the generator set of $\Lambda(R^+)$ from $R^+$ as follows.
If $R^+=\{a,\cdots,b\}$, then we denote the generator set of $\Lambda(R^+)$ by
$\{(a),\cdots,(b)\}$.

The chain complex $(\Lambda(R^+),d)$ is defined as follows.
$d(e)=0$ for all $e\in R^+$.
For $e_1,\cdots,e_n\in R^+$ with $n>1$,
$$\hspace*{-3mm}\begin{array}{l}
\quad d((e_1)(e_2)\cdots(e_n))\vspace{1mm}\\
=\Sigma_{i<j}\,(-1)^{j-i-1}([e_i,e_j])\,
(e_1)\cdots(e_{i-1})(e_{i+1})\cdots
(e_{j-1})(e_{j+1})\cdots(e_n)
\end{array}\quad(3.1)
$$
where $(-e)=-(e)$ and $(0)=0$.
The degree is $|(e_1){\cdots}(e_n)|=n$ and $|1|=0$.

$(\Lambda(R^+),\delta)$ is the dual cochain complex of $(\Lambda(R^+),d)$
regraded as the same group with the dual differential $\delta$.
Then $\delta$ satisfies the Leibnitz formula $\delta(xy)=(\delta x)y+(-1)^{|x|}x(\delta y)$
for all homogeneous $x,y\in\Lambda(R^+)$ and for $e\in R^+$,
$$\delta(e)=\Sigma\,(e')(e''),$$
where the sum is taken over all $e',e''\in R^+$
such that $[e',e'']=e$ (so $\delta(e)=0$ for simple root $e$).

By definition,
$$H_*(R^+)=H_*(\Lambda(R^+),d),\,H^*(R^+)=H^*(\Lambda(R^+),\delta)$$
and for an abelian group $\Bbb F$,\vspace{2mm}\\
\hspace*{13mm}$H_*(R^+;\Bbb F)=H_*(\Lambda(R^+){\otimes}\Bbb F,d),\,H^*(R^+;\Bbb F)
=H^*(\Lambda(R^+){\otimes}\Bbb F,\delta)$.
\vspace{3mm}

{\bf Theorem 3.2}\, {\it $\Lambda(R^+)$ is a weighted (weight graded) group, i.e.,
there is a weight decomposition
$$\Lambda(R^+)=\oplus_{\alpha\in \Omega(R^+)}\,\Lambda(\alpha),$$
where $\Lambda(\alpha)$ is a graded subgroup (called weight subgroup) of $\Lambda(R^+)$.

$(\Lambda(R^+),d)$ and $(\Lambda(R^+),\delta)$ are weighted complexes, i.e., there are complex direct sum
decompositions
$$(\Lambda(R^+),d)=\oplus_{\alpha\in\Omega(R^+)}(\Lambda(\alpha),d),\quad
(\Lambda(R^+),\delta)=\oplus_{\alpha\in\Omega(R^+)}(\Lambda(\alpha),\delta),$$
where $(\Lambda(\alpha),d)$ and $(\Lambda(\alpha),\delta)$
are respectively subcomplex (called weight subcomplex) of $(\Lambda(R^+),d)$ and $(\Lambda(R^+),\delta)$.
The (co)homology group of weight subcomplexes over the abelian group $\Bbb F$ are denoted by
$$H_*^\alpha(R^+;\Bbb F)=H_*(\Lambda(\alpha){\otimes}\Bbb F,d),\,\,
H^*_\alpha(R^+;\Bbb F)=H^*(\Lambda(\alpha){\otimes}\Bbb F,\delta).$$

So $H_*(R^+;\Bbb F)$ and $H^*(R^+;\Bbb F)$ are weighted groups\vspace{2mm}\\
\hspace*{6mm}$H_*(R^+;\Bbb F)=\oplus_{\alpha\in\Omega(R^+)}H_*^\alpha(R^+;\Bbb F),\,\,
H^*(R^+;\Bbb F)=\oplus_{\alpha\in\Omega(R^+)}H^*_\alpha(R^+;\Bbb F)$.
}\vspace{3mm}

{\it Proof}\,\, Let\, $\Bbb Z(\Gamma(R^+))$\, be \,the\, free\, abelian\, group\, generated \, by\, $\Gamma(R^+)$.
Identify $\Bbb Z(\Gamma(R^+))$ with $\Lambda(R^+)$ (only as free abelian groups\,!) as follows.
Suppose a total order is given on $R^+$. Then for
$\sigma\in\Gamma(R^+)\subset\Bbb Z(\Gamma(R^+))$, identify $\sigma$ with the ordered product
$\prod_{e\in \sigma}(e)\in\Lambda(R^+)$.
Specifically, identify $\emptyset\in\Gamma(R^+)\subset\Bbb Z(\Gamma(R^+))$ with $1\in\Lambda(R^+)$.
Then with this identification,
the chain complex $(\Bbb Z(\Gamma(R^+)),d)=(\Lambda(R^+),d)$ is defined.

Precisely, suppose $R^+=\{e_1,\cdots,e_n\}$ with total order $e_1<{\cdots}<e_n$.
For $\sigma\in\Gamma(R^+)$, define the degree $|\sigma|$ to be its cardinality.
For $\{\sigma,\tau\}\in\Gamma(R^+)$ such that $\sigma{\setminus}\tau=\{e_{u},e_{v}\}$,
$\tau{\setminus}\sigma=\{e_{w}\}$, $e_{w}=e_{u}{+}e_{v}$,
define sign function\vspace{2mm}\\
\hspace*{39mm}$\phi(\sigma,\tau)=(-1)^{a+b+c}$, where\\
\hspace*{29mm}$a={\rm the\,\,number\,\,of}\,\,e_k\in\sigma\,\,{\rm such\,\,that}\,\,k\leqslant u,$\\
\hspace*{29mm}$b={\rm the\,\,number\,\,of}\,\,e_k\in\sigma\,\,{\rm such\,\,that}\,\,k\leqslant v,$\hfill(3.2)\\
\hspace*{29mm}$c={\rm the\,\,number\,\,of}\,\,e_k\in\tau\,\,{\rm such\,\,that}\,\,k\leqslant w.$\vspace{1mm}\\
Then for any $\sigma\in\Gamma(R^+)$, we have by (3.1)
$$\hspace{32mm}
\begin{array}{c}
d\sigma=\Sigma_{|\tau|=|\sigma|-1,\,\{\sigma,\tau\}\in\Gamma(R^+)}\,\phi(\sigma,\tau)\tau\vspace{2mm}\\
\delta\sigma=\Sigma_{|\tau|=|\sigma|+1,\,\{\sigma,\tau\}\in\Gamma(R^+)}\,\phi(\sigma,\tau)\tau
\end{array}
\hspace{32mm}(3.3)$$
where the differentials are $0$ if there is no vertex $\tau$ satisfying the condition.

Since for $\{\sigma,\tau\}\in\Gamma(R^+)$, $\omega(\sigma)=\omega(\tau)$,
by (3.3),
$$(\Bbb Z(\Gamma(R^+)),d)=\oplus_{\alpha\in\Omega(R^+)}(\Bbb Z(\Gamma(\alpha)),d),\vspace{-1mm}$$
$$(\Bbb Z(\Gamma(R^+)),\delta)=\oplus_{\alpha\in\Omega(R^+)}(\Bbb Z(\Gamma(\alpha)),\delta).$$

Define $\Lambda(\alpha)=\Bbb Z(\Gamma(\alpha))$ and the theorem holds.
\hfill$\Box$\vspace{3mm}

{\bf Remark} The group action $\circ$ on $\Gamma(R^+)$ makes $\Lambda(R^+)$ a module over
the group ring of $W$.
So there is a group isomorphism $\Lambda(\alpha)\cong
\Lambda(w{*}\alpha)$.
But this group isomorphism is not a (co)chain complex isomorphism, i.e.,
$(\Lambda(\alpha),d)\not\cong(\Lambda(w\alpha),d)$ and
$(\Lambda(\alpha),\delta)\not\cong(\Lambda(w\alpha),\delta)$ in general.
So we never use the group module structure of $\Lambda(R^+)$.
\vspace{3mm}

{\bf Theorem 3.3}\, {\it
Let $\theta\in W$ be the unique element such that $\theta{*}e=-e$ for all $e\in R^+$, equivalently,
$\theta{\circ}\sigma=\sigma^c$ for all $\sigma\in\Gamma(R^+)$.
Then $\theta$ induces a duality complex isomorphism
$$\begin{array}{ccc}
(\Lambda(R^+),d)&\stackrel{\vartheta}{\longrightarrow}&(\Lambda(R^+),\delta)\,\\
\|&&\|\\
\oplus_{\alpha\in\Omega(R^+)}\,(\Lambda(\alpha),d)&
\stackrel{\oplus\vartheta_\alpha}{-\!-\!\!\!\longrightarrow}
&\oplus_{\alpha\in\Omega(R^+)}\,(\Lambda(-\alpha),\delta),
  \end{array}
$$
where $\vartheta_\alpha\colon(\Lambda(\alpha),d)\to(\Lambda(-\alpha),\delta)$
is the restriction of $\vartheta$ on $\Lambda(\alpha)$.
So for any abelian group $\Bbb F$ and weight $\alpha\in\Omega(R^+)$,\vspace{2mm}\\
\hspace*{45mm}$H_*^\alpha(R^+;\Bbb F)=H^*_{-\alpha}(R^+;\Bbb F)$.
}\vspace{3mm}

{\it Proof}\, Let everything be as in the proof of Theorem 3.2
and $(\Bbb Z(\Gamma(R^+)),\delta)$ be the dual complex of $(\Bbb Z(\Gamma(R^+)),d)$.

For $\sigma=\{e_{i_1},\cdots,e_{i_s}\}\in\Gamma(R^+)$,
define $\phi(\sigma)=(-1)^{i_1+\cdots+i_s}$ and $\phi(\emptyset)=1$.
We prove that for any edge $\{\sigma,\tau\}\in\Gamma(R^+)$,
$$\phi(\sigma)\phi(\tau)\phi(\sigma,\tau)=\phi(\sigma^c,\tau^c).$$
Suppose $\sigma{\setminus}\tau=\{e_{u},e_{v}\}$, $\tau{\setminus}\sigma=\{e_{w}\}$.
Then $\phi(\sigma,\tau)\phi(\sigma^c,\tau^c)=(-1)^{u{+}v{+}w}$ and
$\phi(\sigma)\phi(\tau)=\phi(\sigma^c)\phi(\tau^c)=(-1)^{u+v+w}$. So the equality holds.

For any $\sigma\in\Gamma(R^+)\subset\Bbb Z(\Gamma(R^+))$, define
$$\vartheta(\sigma)=\phi(\sigma)\,(\theta{\circ}\sigma)=\phi(\sigma)\,\sigma^c.$$
Then
\begin{eqnarray*}&&\vartheta(d\sigma)\\
&=&\vartheta\big(\,\Sigma_{|\tau|=|\sigma|-1,\,\{\sigma,\tau\}\in\Gamma(R^+)}\,\,\phi(\sigma,\tau)\tau\,\big)\\
&=&\Sigma_{|\sigma|=|\tau|+1,\,\{\sigma,\tau\}\in\Gamma(R^+)}\phi(\tau)\phi(\sigma,\tau)\tau^c\\
&=&\Sigma_{|\sigma|=|\tau|+1,\,\{\sigma,\tau\}\in\Gamma(R^+)}\phi(\sigma)\phi(\sigma^c,\tau^c)\tau^c\\
&=&\phi(\sigma)\delta \sigma^c\\
&=&\delta\vartheta(\sigma)
\end{eqnarray*}
and $\vartheta$ is a dual complex isomorphism.
Since $\omega(\theta{\circ}\sigma)=\omega(\sigma^c)=-\omega(\sigma)$, $\vartheta_\alpha$
is an isomorphism from $(\Lambda(\alpha),d)$ to $(\Lambda(-\alpha),\delta)$.
\hfill $\Box$\vspace{3mm}

{\bf Theorem 3.4} {\it For $\alpha\in\Omega(R^+)$, $\Lambda(\alpha)\cong\Bbb Z$ if and only if $\alpha$ satisfies the equivalent conditions in Theorem 2.5.
}\vspace{3mm}

{\it Proof}\, $\Lambda(\alpha)\cong\Bbb Z$ if and only if $\Gamma(\alpha)$ has only one vertex.
\hfill$\Box$\vspace{3mm}

\section{Index and neighbors}\vspace{3mm}

\hspace*{5mm}
{\bf Definition 4.1} For a weight subgraph $\Gamma(\alpha)$ and $e\in R^+$,
$\Gamma_e^1(\alpha)$ is the full subgraph of $\Gamma(\alpha)$ consisting of all vertex $\sigma$
such that $e\in\sigma$ and
$\Gamma_e^0(\alpha)$ is the full subgraph of $\Gamma(\alpha)$ consisting of all vertex $\sigma$ such that $e\not\in \sigma$.
\vspace{3mm}

{\bf Theorem 4.2} {\it For any $e\in R^+$ and $\alpha\in\Omega(R^+)$,
there is an integer $r(e,\alpha)\in\Bbb Z$ such that
$r_e{*}\alpha-\alpha=r(e,\alpha)e$. $r(e,\alpha)$ is called the index of $\alpha$
with respect to $e$.
}\vspace{3mm}

{\it Proof}\, Take any $\sigma\in\Gamma(\alpha)$, then $\Sigma_{f\in\sigma}\,(r_e{*}f-f)=ke$ with $k$ an integer.
Since $r_e{*}\varrho-\varrho=-e$, we have
$$\hspace{12mm}r{*}\alpha-\alpha=r_e{*}\varrho-(\Sigma_{f\in\sigma}\,r_e{*}f)
-\big(\varrho-(\Sigma_{f\in\sigma}\,f)\big)=-(k{+}1)e
\hspace{12mm}(4.1)$$
is an integer.
\hfill$\Box$\vspace{3mm}

{\bf Lemma 4.3} {\it
For a simple root $e$ and any weight $\alpha\in\Omega(R^+)$, we have the following.\vspace{1mm}

(1) If $r(e,\alpha)\geqslant 0$, then $\Gamma_e^1(\alpha)\neq\emptyset$.
Moreover, if $\Gamma_e^0\neq\emptyset$, then
for any vertex $\sigma\in\Gamma_e^0(\alpha)$, there is $\tau\in\Gamma_e^1(\alpha)$ such that $\{\sigma,\tau\}\in\Gamma(\alpha)$.
\vspace{1mm}

(2) If $r(e,\alpha)\leqslant 0$, then $\Gamma_e^0(\alpha)\neq\emptyset$.
Moreover, if $\Gamma_e^1\neq\emptyset$, then
for any vertex $\sigma\in\Gamma_e^1(\alpha)$, there is $\tau\in\Gamma_e^0(\alpha)$ such that $\{\sigma,\tau\}\in\Gamma(\alpha)$.
}\vspace{3mm}

{\it Proof}\,
The $*$ action of $r_e$ on $R^+{\setminus}\{e\}$ is a permutation and
so $R^+{\setminus}\{e\}$ is in one of the following cases.

i) $R^+{\setminus}\{e\}=\{e_1,{\cdots},e_s,\o e_1,{\cdots},\o e_s,h_1,{\cdots},h_t\}$,
where $r_e{*}e_i=\o e_i=e_i+e$ and $r_e{*}h_i=h_i$.
In this case, $e$ is the longer root in the irreducible component if there are roots of different lengths
in the component.

ii) $R^+{\setminus}\{e\}=\{e_1,\!{\cdots},e_s,\o e_1,\!{\cdots},\o e_s,f_1,\!{\cdots},f_s,\o f_1,\!{\cdots},\o f_s,
f'_1,\!{\cdots},f'_s,h_1,\!{\cdots},h_t\}$,
where $r_e{*}e_i=\o e_i=e_i{+}e$, $r_e{*}f_i=\o f_i=f_i{+}2e$, $r_e{*}f'_i=f'_i$, $f'_i=f_i{+}e$ and $r_e{*}h_i=h_i$.
In this case, $e$ is the shorter root in the irreducible component with
Dynkin graph other than $G_2$.

iii) $R^+{\setminus}\{e\}=\{e_1,\o e_1,g_1,\o g_1,h_1,{\cdots},h_s\}$,
where $r_e{*}e_1=\o e_1=e_1{+}e$, $r_e{*}g_1=\o g_1=g_1{+}3e$, $\o g_1=\o e_1{+}e$ and $r_e{*}h_i=h_i$.
In this case, $e$ is the shorter root in the irreducible component with Dynkin graph $G_2$.

We prove conclusion (1) of the theorem in the following.\vspace{2mm}

Case i) Suppose $\Gamma_e^1(\alpha)=\emptyset$ and $\sigma\in\Gamma_e^0(\alpha)$.
Then $\o e_i\in \sigma$ implies $e_i\in \sigma$, for otherwise $(\sigma{\setminus}\{\o e_i\}){\cup}\{e_i,e\}\in\Gamma_e^1(\alpha)$.
So
$$\sigma=\{e_{i_1},{\cdots},e_{i_u}, \o e_{i_1},{\cdots},\o e_{i_u}, e_{j_1},\cdots,e_{j_v},
h_{k_1},\cdots,h_{k_w}\}$$
and by (4.1), $r_e{*}\alpha{-}\alpha=-(1{+}v)e$.
A contradiction\,! So $\Gamma_e^1(\alpha)\neq\emptyset$.

Suppose  $\Gamma_e^0(\alpha)\neq \emptyset$, then for $\sigma\in\Gamma_e^0(\alpha)$,
$$\sigma=\{e_{i_1},{\cdots},e_{i_u}, \o e_{j_1},{\cdots},\o e_{j_v},h_{k_1},
\cdots,h_{k_w}\}.$$
By (4.1), $r_e{*}\alpha{-}\alpha=-(1{+}u{-}v)e$. So $v>u$ and there is $\o e_{j_s}\in \sigma$ such that
$e_{j_s}\not\in \sigma$.
So $\tau=(\sigma{\setminus}\{\o e_{j_s}\})\cup\{e_{j_s},e\}\in\Gamma_e^1(\alpha)$.\vspace{2mm}

Case ii) Suppose $\Gamma_e^1(\alpha)=\emptyset$ and $\sigma\in\Gamma_e^0(\alpha)$.
Then for the same reason as in case i),
$\o e_i\in \sigma$ implies $e_i\in \sigma$,
$f'_i\in \sigma$ implies $f_i\in \sigma$ and $\o f_i\in \sigma$ implies $f'_i\in \sigma$.
So\vspace{2mm}\\
\hspace*{5mm}$\sigma=\{e_{i_1},{\cdots},e_{i_u}$, $\o e_{i_1},{\cdots},\o e_{i_u}$,
$f_{j_1},{\cdots},f_{j_v}$, $f'_{j_1},{\cdots},f'_{j_v}$,
$f'_{k_1},{\cdots},f'_{k_w}$,\\
\hspace*{35mm}$\o f_{k_1},{\cdots},\o f_{k_w}$,
$e_{l_1},{\cdots},e_{l_x}$,
$f_{m_1},{\cdots},f_{m_y}$,
$h_{n_1},\cdots,h_{n_z}\}$,\vspace{2mm}\\
where $\{k_1,{\cdots},k_w\}\subset\{j_1,{\cdots},j_v\}$.
So $v>w$ and by (4.1)\vspace{1mm}\\
\hspace*{38.5mm}$r_e{*}\alpha{-}\alpha=-(1{+}x{+}2y{+}2(v{-}w))e.$\vspace{1mm}\\
A contradiction\,!
So $\Gamma_e^1(\alpha)\neq\emptyset$.

Suppose $\Gamma_e^0(\alpha)\neq\emptyset$ and $\sigma\in\Gamma_e^0(\alpha)$.
Then\vspace{2mm}\\
\hspace*{5mm}$\sigma=\{e_{i_1},{\cdots},e_{i_u}$, $\o e_{j_1},{\cdots},\o e_{j_v}$,
$f_{k_w},{\cdots},f_{k_w}$, $f'_{l_1},{\cdots},f'_{l_x}$,\\
\hspace*{85mm}$\o f_{m_1},{\cdots},\o f_{m_y}$,
$h_{n_1},\cdots,h_{n_z}\}$,\vspace{2mm}\\
By (4.1), $r_e{*}\alpha{-}\alpha=-(1{+}u{-}v{+}2(w{-}y))e$.
If $u>v$, then $y>w$. So either there is $\o e_{i_s}\in \sigma$ such that $e_{i_s}\not\in \sigma$,
or there is $\o f_{k_t}\in \sigma$ such that $f_{k_t}\not\in \sigma$.
If $\o e_{i_s}\in \sigma$ but $e_{i_s}\not\in \sigma$,
then $\tau=(\sigma{\setminus}\{\o e_{i_s}\})\cup\{e_{i_s},e\}\in\Gamma_e^1(\alpha)$.
If $\o f_{k_t}\in \sigma$ but $f_{k_t}\not\in \sigma$, then $\tau=(\sigma{\setminus}\{\o f_{k_t}\})\cup\{f_{k_t},e\}\in\Gamma_e^1(\alpha)$.\vspace{2mm}

Case iii) Suppose $\Gamma_e^1(\alpha)=\emptyset$ and $\sigma\in\Gamma_e^0(\alpha)$.
For the same reason as in case i), $\o e_1\in \sigma$ implies $e_1\in \sigma$, $\o g_1\in \sigma$ implies $\o e_1\in \sigma$,
$e_1\in \sigma$ implies $g_1\in \sigma$. So $\sigma$ is one of the following (where the omitted part consists of $h_i$)
$$\{e_1,\cdots\},\{g_1,\cdots\},\{e_1,g_1,\cdots\},\{e_1,\o e_1,g_1,\cdots\},\{e_1,\o e_1,g_1,\o g_1,\cdots\}.$$
In all the above cases, $r(e,\alpha)<0$.
A contradiction! So $\Gamma_e^1(\alpha)\neq\emptyset$.

Suppose $\Gamma_e^0(\alpha)\neq\emptyset$ and $\sigma\in\Gamma_e^0(\alpha)$.
Then $\sigma$ is one of the following (where the omitted part consists of $h_i$)
$$\{\o e_1,\cdots\},\{\o g_1,\cdots\},\{e_1,\o g_1,\cdots\},\{\o e_1,g_1,\o g_1,\cdots\},\{e_1,\o e_1,\o g_1,\cdots\}.$$
The corresponding $\tau\in\Gamma_e^1(\alpha)$ such that $\{\sigma,\tau\}\in\Gamma(\alpha)$ is
$$\{e,e_1,\cdots\},\{e,\o e_1,\cdots\},\{e,e_1,\o e_1,\cdots\},\{e,e_1,g_1,\o g_1,\cdots\},\{e,\o e_1,g_1,\o g_1,\cdots\}.$$

Conclusion (2) of the theorem follows by duality isomorphism $\sigma\to \theta{\circ}\sigma$ in Theorem 3.3.
If $r(e,\alpha)\leqslant 0$,
then $r(e,-\alpha)\geqslant 0$ and conclusion (1) holds for $\Gamma({-}\alpha)$.
So conclusion (2) holds for $\Gamma(\alpha)$ by the 1-1 correspondence.
\hfill$\Box$\vspace{3mm}

{\bf Theorem 4.4} {\it Lemma 4.3 holds for all roots, i.e.,
for any $e\in R^+$ and weight $\alpha\in\Omega(R^+)$, we have the following.\vspace{1mm}

(1) If $r(e,\alpha)\geqslant 0$, then $\Gamma_e^1(\alpha)\neq\emptyset$.
Moreover, if $\Gamma_e^0\neq\emptyset$, then
for any vertex $\sigma\in\Gamma_e^0(\alpha)$, there is $\tau\in\Gamma_e^1(\alpha)$ such that $\{\sigma,\tau\}\in\Gamma(\alpha)$.
\vspace{1mm}

(2) If $r(e,\alpha)\leqslant 0$, then $\Gamma_e^0(\alpha)\neq\emptyset$.
Moreover, if $\Gamma_e^1\neq\emptyset$, then
for any vertex $\sigma\in\Gamma_e^1(\alpha)$, there is $\tau\in\Gamma_e^0(\alpha)$ such that $\{\sigma,\tau\}\in\Gamma(\alpha)$.
}\vspace{3mm}

{\it Proof}\,
Suppose $e$ is not a simple root.
Then there is $w\in W$ such that $w{*}e$ is a simple root.
From the equality $r_e{*}\alpha{-}\alpha=r(e,\alpha)e$, we have
$$r_{w*e}{*}(w{*}\alpha){-}w{*}\alpha=(wr_ew^{-1}){*}(w{*}\alpha){-}w{*}\alpha
=w{*}(r_e{*}\alpha{-}\alpha)=r(e,\alpha)w{*}e.$$
So $r(w{*}e,w{*}\alpha)=r(e,\alpha)$. Since $\Gamma(w{*}\alpha)\cong\Gamma(\alpha)$ and $\Gamma_{w{*}e}^i(w{*}\alpha)\cong\Gamma_e^i(\alpha)$,
that the theorem holds for $w{*}e$ implies that the theorem holds for $e$.
\hfill$\Box$\vspace{3mm}

\section{Basic properties of basis graph}\vspace{3mm}

\hspace*{5mm}
{\bf Definition 5.1} A diamond graph $\Gamma$ is a graph satisfying the following condition.
If there are $a,b,c\in\Gamma$ such that $\{a,b\}\in\Gamma$ and $\{b,c\}\in\Gamma$,
then there is a unique $d\in\Gamma$, $d\neq b$ such that $\{a,d\}\in\Gamma$, $\{d,c\}\in\Gamma$ and
$\{a,c\}\not\in\Gamma$, $\{b,d\}\not\in\Gamma$.
The full subgraph with vertex set $\{a,b,c,d\}$ is called a diamond of $\Gamma$.
We denote by $\langle a;b;c;d\rangle$ the diamond such that
$\{a,b\},\{b,c\},\{c,d\},\{d,a\}\in\Gamma$ and $\{a,c\},\{b,d\}\not\in\Gamma$.\vspace{3mm}

By the above definition, a diamond graph has no triangular subgraph.\vspace{3mm}

{\bf Theorem 5.2}\, {\it $\Gamma(R^+)$ is a diamond graph.
}\vspace{3mm}

{\it Proof}\, Let $\sigma,\tau,\mu\in \Gamma(\alpha)$ satisfy $\{\sigma,\tau\},\{\tau,\mu\}\in\Gamma(\alpha)$.
We may suppose $\sigma,\tau,\mu$ is one of the following two cases.

i) $|\sigma|{+}2=|\tau|{+}1=|\mu|$.

ii) $|\sigma|=|\mu|=|\tau|{\pm}1$.

Case i)
Suppose $\sigma{\setminus}\tau=\{e\}$, $\tau{\setminus}\sigma=\{e_1,e_2\}$ with $e=e_1{+}e_2$,
$\tau{\setminus}\mu=\{f\}$, $\mu{\setminus}\tau=\{f_1,f_2\}$ with $f=f_1{+}f_2$.
We have the following two cases.

(1) $f\not\in\{e_1,e_2\}$.
Then we have a diamond $\langle \sigma;\tau;\mu;\nu\rangle$
as in the following figure
(where the arrow means an edge pointing to the vertex with bigger degree
and all the omitted parts are the same).
\[\begin{array}{c}
\sigma=\{e,f,\cdots\}\\
\swarrow\hspace{25mm}\searrow\\
\tau=\{e_1,e_2,f,\cdots\}\quad \nu=\{e,f_1,f_2,\cdots\}\\
\searrow\hspace{25mm}\swarrow\\
\mu=\{e_1,e_2,f_1,f_2,\cdots\}
  \end{array}
\]

Now we prove that $\nu$ is unique.
Suppose there is a vertex $\nu'\neq \nu$ such that $\langle \sigma;\tau;\mu;\nu'\rangle$ is a diamond.
Let $\sigma{\setminus}\nu'=\{g\}$, $\nu'{\setminus}\sigma=\{g_1,g_2\}$ with $g=g_1{+}g_2$.
If $e=g$, then $\nu'=\{f,g_1,g_2,\cdots\}$ with $\{g_1,g_2\}{\cap}\{e_1,e_2\}=\emptyset$.
It is obvious $\{\nu',\mu\}\not\in\Gamma(R^+)$\,!
So $e\neq g$. Analogously, we have $f\neq g$.
If $g\neq e$ and $g\neq f$,
then $\nu'=\{e,f,{\cdots}\}$ and obviously $\{\nu',\mu\}\not\in\Gamma(R^+)$\,!
So $\nu$ is unique.

(2) $f\in\{e_1,e_2\}$. We may suppose $f=e_2$.
Since
$$[e_1,[f_1,f_2]]+[f_1,[f_2,e_1]]+[f_2,[e_1,f_3]]=0$$
and $e_1{+}f_1{+}f_2\in R^+$, we have that one of $e_1{+}f_1,e_1{+}f_2$
is a positive root and the other is not a root.
We may suppose $e_1{+}f_2$ is not a root.
Then we have a diamond $\langle \sigma;\tau;\mu;\nu\rangle$
as in the following figure.
\[\begin{array}{c}
\sigma=\{e,\cdots\}\\
\swarrow\hspace{25mm}\searrow\\
\tau=\{e_1,f_1{+}f_2,\cdots\}\quad \nu=\{e_1{+}f_1,f_2,\cdots\}\\
\searrow\hspace{25mm}\swarrow\\
\mu=\{e_1,f_1,f_2,\cdots\}
  \end{array}
\]

Now we prove $\nu$ is unique.
Suppose there is $\nu'\neq \nu$ such that $\langle \sigma;\tau;\mu;\nu'\rangle$ is a diamond.
Let $\sigma{\setminus}\nu'=\{h\}$, $\nu'{\setminus}\sigma=\{h_1,h_2\}$ with $h=h_1{+}h_2$.
If $h=e$, then $\nu=\{h_1,h_2,\cdots\}$ and so
$\{h_1,h_2\}{\cap}\{e_1,f_1{+}f_2,e_1{+}f_1,f_2\}=\emptyset$.
$\{\nu',\mu\}\in\Gamma(R^+)$ implies that at least one of $h_1$ and $h_2$ is in $\{e_1,f_1,f_2\}$.
We may suppose $h_1\in\{e_1,f_1,f_2\}$.
Since $\{h_1,h_2\}{\cap}\{e_1,f_1{+}f_2,e_1{+}f_1,f_2\}=\emptyset$, we have $h_1=f_1$.
Then from $h_1{+}h_2=e_1{+}f_1{+}f_2$ we have $h_2=e_1{+}f_2$.
But $e_1{+}f_2$ is not a root. A contradiction!
So $h\neq e$.
If $h\neq e$, then $e\in \nu'$.
$\{\nu',\mu\}\in\Gamma(R^+)$ implies that either $e\in \mu$, or $\{e_1,e_2\}\in \mu$. A contradiction!
So $\nu$ is unique.

Case ii) We prove that there is $w\in W$ such that $|w{\circ}\sigma|,|w{\circ}\tau|,|w{\circ} \mu|$ is in case i).
If so, there is a unique $\nu'$ such that $\langle w{\circ}\sigma,w{\circ}\tau,w{\circ}\mu,\nu'\rangle$ is a diamond.
So $\langle \sigma,\tau,\mu,w^{-1}{\circ}\nu'\rangle$ is a diamond.

We first prove a lemma. Suppose for simple reflection $r_e$ and $\{\sigma,\tau\}\in\Gamma(R^+)$,
we have $|\sigma|=|\tau|$ and $|r_e{\circ}\sigma|=|r_e{\circ}\tau|$.
Then $f\in\sigma$ and $f\not\in\tau$ implies $r_e{*}f\in r_e{\circ}\sigma$ and $r_e{*}f\not\in r_e{\circ}\tau$.
$|\sigma|=|\tau|$ and $|r_e{\circ}\sigma|=|r_e{\circ}\tau|$ imply either $e\in\sigma{\cap}\tau$,
or $e\in\sigma^c{\cap}\tau^c$. So $f\in\sigma$ and $f\not\in\tau$ imply $f\neq e$.
So if $r_e{*}f\in r_e{\circ}\tau$, then $f\in\tau$. A contradiction!
The lemma is proved.

Suppose for all $w'\in W$, $|w'{\circ}\sigma|=|w'{\circ}\mu|$.
Take $f\in\sigma$ but $f\not\in\tau$.
There is $w\in W$ such that $s=w{*}e$ is a simple root.
Suppose $w=r_{e_1}\cdots r_{e_n}$, where $e_1,{\cdots},e_n$ (repetition allowed) are all simple root.
Denote
$$\sigma_0=\sigma,\,\,\tau_0=\tau,\,\,\mu_0=\mu,\,\,f_0=f,\vspace{-1mm}$$
$$\sigma_i=r_{e_i}{\circ}\sigma_{i-1},\,\,\tau_i=r_{e_i}{\circ}\tau_{i-1},\,\,
\mu_i=r_{e_i}{\circ}\mu_{i-1},\,\,f_i=r_{e_i}{*}f_{i-1},$$
for $i=1,\cdots,n$.
By the lemma and induction hypothesis $|\sigma_i|=|\mu_i|$,
$f_i\in\sigma_i$ but $f_i\not\in\mu_i$ for all $i$.
Specifically, $s=f_n\in w{\circ}\sigma$ but $s\not\in w{\circ}\mu$.
Then $|r_s{\circ}(w{\circ}\sigma)|=|w{\circ}\sigma|{-}1$ and $|r_s{\circ}(w{\circ}\mu)|=|w{\circ}\mu|{+}1$.
$|(r_sw){\circ}\sigma|\neq|(r_sw){\circ}\mu|$. A contradiction!
So there must be $w\in W$ such that $|w{\circ}\sigma|\neq|w{\circ}\mu|$, equivalently,
$|w{\circ}\sigma|=|w{\circ}\mu|{\pm}2$.
So $|w{\circ}\sigma|,|w{\circ}\tau|,|w{\circ} \mu|$ is in case i).
\hfill$\Box$\vspace{3mm}

{\bf Theorem 5.3}\, {\it The weight decomposition
$$\Gamma(R^+)=\sqcup_{\alpha\in\Omega(R^+)}\,\Gamma(\alpha)$$
in Theorem 2.4 is the connected component decomposition, i.e., every
weight subgraph $\Gamma(\alpha)$ is a connected graph.
}\vspace{3mm}

{\it Proof}\,
Let $\Omega=\{w\varrho\,|\,w\in W\}\subset\Omega(R^+)$.
By Theorem 2.5, the theorem holds for $\alpha\in\Omega$.

Suppose $\alpha\not\in\Omega$, then by 4$.$ of Theorem 2.5, there is $e_1\in R^+$ such that $r_{e_1}{*}\alpha=\alpha$.
Define $\alpha_1=\alpha{-}e_1$. Inductively define
$\alpha_{i}=\alpha_{i-1}-e_i$ ($\alpha_0=\alpha$) if $\alpha_{i-1}\not\in\Omega$
and $r_{e_i}{*}\alpha_{i-1}=\alpha_{i-1}$.
Since $\alpha_i=\alpha{-}e_1{-}{\cdots}{-}e_{i}$
and $\Omega(R^+){\setminus}\Omega$ is finite, there is an $n$ such that $\alpha_n\in\Omega$
and $\alpha,\alpha_1,{\cdots},\alpha_{n-1}\not\in\Omega$.

Suppose the following lemma holds. There are $\beta_1,{\cdots},\beta_{n-1}\not\in\Omega$, $\beta_n\in\Omega$,
and positive roots $f_1,{\cdots},f_n$ satisfying the following property. $\beta_1=\alpha_1$ and $f_1=e_1$.
For $i\geqslant 1$, $r_{f_{i+1}}{*}\beta_i=\beta_i$,
$$\beta_{i+1}=\left\{\begin{array}{ccc}
\beta_i-f_{i+1}&&{\rm if}\,f_{i+1}\not\in F_i,\\
\beta_i+f_{i+1}&&{\rm if}\,f_{i+1}\in F_i,
\end{array}\right.$$
where $F_1=\{f_1\}$, $F_{i+1}=\left\{
\begin{array}{ll}
F_i{\cup}\{f_{i+1}\}&{\rm if}\,f_{i+1}\not\in F_i,\vspace{1mm}\\
F_i\,{\setminus}\{f_{i+1}\}&{\rm if}\,f_{i+1}\in F_i.
\end{array}\right.$\vspace{1mm}

Then we prove that for any vertex $\sigma\in\Gamma(\alpha)$,
there are $\sigma_1,{\cdots},\sigma_n\in\Gamma(\alpha)$
such that either $\sigma_{i-1}=\sigma_i$, or $\{\sigma_{i-1},\sigma_{i}\}\in\Gamma(\alpha)$, $F_i\subset \sigma_i$ 
and $\sigma_i{\setminus}F_i\in\Gamma(\beta_i)$ for $i=1,{\cdots},n$ ($\sigma_0=\sigma$).

If $\sigma\in\Gamma_{e_1}^1(\alpha)$, then define $\sigma_1=\sigma$.
If $\sigma\in\Gamma_{e_1}^0(\alpha)$, then since $r_{f_1}{*}\alpha=\alpha$, by Theorem 4.4,
there is $\sigma_1\in\Gamma_{f_1}^1(\alpha)$ such that $\{\sigma,\sigma_1\}\in\Gamma(\alpha)$.
So $\sigma_1\in\Gamma(\alpha)$ is defined such that $F_1\subset \sigma_1$ and $\sigma_1{\setminus}F_1\in\Gamma(\beta_1)$.

Suppose $\sigma_i\in\Gamma(\alpha)$ is defined. Then define $\sigma_{i+1}$ as follows.\vspace{1mm}

(1) $f_{i+1}\not\in F_i$ and $f_{i+1}\in \sigma_i$. Define $\sigma_{i+1}=\sigma_i$.
Then $F_{i+1}\subset \sigma_{i+1}$ and $\sigma_{i+1}{\setminus}F_{i+1}\in\Gamma(\beta_{i}{-}f_{i+1})
=\Gamma(\beta_{i+1})$.\vspace{1mm}

(2) $f_{i+1}\in F_i$ and so $f_{i+1}\in \sigma_i$. Define $\sigma_{i+1}=\sigma_i$.
Then $F_{i+1}\subset \sigma_{i+1}$ and $\sigma_{i+1}{\setminus}F_{i+1}
=\sigma_i{\setminus}F_{i+1}\in\Gamma(\beta_i{+}f_{i+1})=\Gamma(\beta_{i+1})$.\vspace{1mm}

(3) $f_{i+1}\not\in F_i$ and $f_{i+1}\not\in \sigma_i$.
Since $r_{f_{i+1}}{*}\beta_i=\beta_i$, by Theorem 4.4,
there is $\sigma'_{i+1}\in\Gamma_{f_{i+1}}^1(\beta_i)$ with 
$\{\sigma_i{\setminus} F_i,\sigma'_{i+1}\}\in\Gamma(\beta_i)$.
Define $\sigma_{i+1}=F_i{\cup}\sigma'_{i+1}$.
Then $\{\sigma_i,\sigma_{i+1}\}\in\Gamma(\alpha)$, $F_{i+1}\subset \sigma_{i+1}$ and
$\sigma_{i+1}{\setminus}F_{i+1}=\sigma'_{i+1}{\setminus}\{f_{i+1}\}\in\Gamma(\beta_{i}{-}f_{i+1})
=\Gamma(\beta_{i+1})$.\vspace{1mm}

Thus, for any $\sigma\in\Gamma(\alpha)$, $\sigma_1,{\cdots},\sigma_n$ are defined.
For any $\sigma,\tau\in\Gamma(\alpha)$,
we have a sequence $\sigma,\sigma_1,{\cdots},\sigma_{n-1}, \sigma_n,\tau_n,\tau_{n-1},{\cdots},\tau_1,\tau$ (repetition allowed).
Since $\omega(\sigma_n{\setminus}F_n)=\omega(\tau_n{\setminus}F_n)=\beta_n\in\Omega$,
we have $\sigma_n{\setminus}F_n=\tau_n{\setminus}F_n$ and $\sigma_n=\tau_n$.
So the above sequence is a path
connecting $\sigma$ and $\tau$. $\Gamma(\alpha)$ is connected.

The lemma is proved by inductively define $\beta_i\!=\!w_i{*}\alpha_i$ with all $w_i\in W$ as follows.
$w_1=1$.
$$f_{i+1}=\left\{
\begin{array}{ll}
\hspace*{3mm} w_{i}{*}e_{i+1}&{\rm if}\,w_{i}{*}e_{i+1}\in R^+,\vspace{1mm}\\
{-}w_{i}{*}e_{i+1}&{\rm if}\,w_{i}{*}e_{i+1}\in R^-,
\end{array}\right.$$
$$\hspace*{20mm}w_{i+1}=\left\{
\begin{array}{ll}
w_i&{\rm if}\,\,f_{i+1}\not\in F_i\,\,{\rm and}\,\,w_i{*}e_{i+1}\in R^+,\vspace{1mm}\\
w_ir_{e_{i+1}}&{\rm if}\,\,f_{i+1}\not\in F_i\,\,{\rm and}\,\,w_i{*}e_{i+1}\in R^-,\vspace{1mm}\\
w_ir_{e_{i+1}}&{\rm if}\,\,f_{i+1}\in F_i\,\,{\rm and}\,\,w_i{*}e_{i+1}\in R^+,\vspace{1mm}\\
w_i&{\rm if}\,\,f_{i+1}\in F_i\,\,{\rm and}\,\,w_i{*}e_{i+1}\in R^-.\vspace{1mm}
\end{array}\right.$$

If $f_{i+1}\not\in F_i$, then
$$\beta_i{-}f_{i+1}=\left\{
\begin{array}{ll}
w_i{*}\alpha_i{-}w_i{*}e_{i+1}=w_{i}{*}\alpha_{i+1}=\beta_{i+1}&{\rm if}\,\, w_i{*}e_{i+1}\in R^+,\vspace{1mm}\\
w_i{*}\alpha_i{+}w_i{*}e_{i+1}=w_ir_{e_{i+1}}{*}\alpha_{i+1}=\beta_{i+1}&{\rm if}\,\,w_i{*}e_{i+1}\in R^-.
\end{array}\right.
$$

If $f_{i+1}\in F_i$, then
$$\beta_i{+}f_{i+1}=\left\{
\begin{array}{ll}
w_i{*}\alpha_i{+}w_i{*}e_{i+1}=w_{i}r_{e_{i+1}}{*}\alpha_{i+1}=\beta_{i+1}&{\rm if}\,\, w_i{*}e_{i+1}\in R^+,\vspace{1mm}\\
w_i{*}\alpha_i{-}w_i{*}e_{i+1}=w_i{*}\alpha_{i+1}=\beta_{i+1}&{\rm if}\,\,w_i{*}e_{i+1}\in R^-.
\end{array}\right.
$$

$r_{f_{i+1}}{*}\beta_i=(w_ir_{e_{i+1}}w^{-1}_i){*}\beta_i=(w_ir_{e_{i+1}}w^{-1}_i){*}(w_i{*}\alpha_i)=w_i{*}\alpha_i
=\beta_i$.\vspace{1mm}
\hfill$\Box$\vspace{3mm}

{\bf Definition 5.4} A diamond graph $\Gamma$ is admissible if there is a sign function
$\phi$ defined from all edges of $\Gamma$ to the set $\{\pm 1\}$ such that for every diamond $\langle a;b;c;d\rangle$,
$\phi(a,b)\phi(b,c)+\phi(a,d)\phi(d,c)=0$.
\vspace{3mm}

{\bf Theorem 5.5} {\it $\Gamma(R^+)$ is an admissible diamond graph with sign function
$\phi$ as defined in {\rm (3.2)}.
}\vspace{3mm}

{\it Proof}\,
Notice that the sign function $\phi$ depends on the total order on $R^+$.
We have to consider sign functions defined with respect to different orders on $R^+$.
Suppose $R^+=\{e_1,\cdots,e_n\}$ with total order $e_1<{\cdots}<e_n$ and denote by $\phi$ the sign function
defined with respect to this total order.
For a permutation $\pi$ on $\{1,{\cdots},n\}$, there is a new total order $e_{\pi(1)}<{\cdots}<e_{\pi(n)}$
on $R^+$.
Denote by $\phi_\pi$ the sign function defined with respect to this total order.

There are two types of diamonds in $\Gamma(R^+)$.
A diamond $\langle \sigma;\tau;\mu;\nu\rangle$ is called flat if $|\sigma|{-}|\mu|\neq |\tau|{-}|\nu|$.
A diamond $\langle \sigma;\tau;\mu;\nu\rangle$ is called folded if $|\sigma|{-}|\mu|=|\tau|{-}|\nu|=0$.
For a flat diamond $\langle \sigma;\tau;\mu;\nu\rangle$ such that $|\sigma|{+}2=|\tau|{+}1=|\mu|$,
by (3.3),
$$\delta^2 \sigma=
\Sigma_{i,\,j,\,|\sigma''_j|=|\sigma|+2}\,\phi(\sigma,\sigma'_i)\phi(\sigma'_i,\sigma''_j)\sigma''_j=0$$
we have that $\phi(\sigma,\tau)\phi(\tau,\mu)+\phi(\sigma,\nu)\phi(\nu,\mu)=0$ for some $i,j,k$ such that
$\tau=\sigma'_i$, $\nu=\sigma'_j$ and $\mu=\sigma''_k$.
For the same reason, for any permutation $\pi$ on $\{1,{\cdots},n\}$ and flat diamond
$\langle \sigma;\tau;\mu;\nu\rangle$,
$$\hspace{28mm}
\phi_\pi(\sigma,\tau)\phi_\pi(\tau,\mu)+\phi_\pi(\sigma,\nu)\phi_\pi(\nu,\mu)=0.\hspace{28mm}(5.1)$$

For a simple root $s$, let $\pi_s$ be the permutation on $\{1,{\cdots},n\}$ defined by
$\pi_s(k)=k$ if $s=e_k$ and $\pi_s(i)=j$ if $r_s{*}e_i=e_j$ for $i\neq k$.
Beware that these permutations $\pi_s$ do not induce a group action of $W$ on $\{1,{\cdots},n\}$.
By (3.2), for any $\{\xi,\eta\}\in\Gamma(R^+)$ and any permutation $\pi$ on
$\{1,{\cdots},n\}$,
$$\hspace{38mm}\phi_\pi(r_s{\circ}\xi,r_s{\circ}\eta)=\phi_{\pi\pi_s}(\xi,\eta).\hspace{38mm}(5.2)$$

Suppose $\langle \sigma;\tau;\mu;\nu\rangle$ is a folded diamond such that $|\sigma|=|\mu|=|\tau|{\pm}1$.
From the proof of case ii) of Theorem 5.2, there is $w\in W$ such that
$\langle w{\circ}\sigma;w{\circ}\tau;w{\circ}\mu;w{\circ}\nu\rangle$ is a flat diamond.
Suppose $w=r_{s_1}{\cdots}r_{s_t}$,
where $s_1,{\cdots},s_t$ are simple roots. Let $\pi=\pi_{s_t}{\cdots}\pi_{s_1}$.
Then by (5.1)
$$\phi_\pi(w{\circ}\sigma,w{\circ}\tau)\phi_\pi(w{\circ}\tau,w{\circ}\mu)+
\phi_\pi(w{\circ}\sigma,w{\circ}\nu)\phi_\pi(w{\circ}\nu,w{\circ}\mu)=0.$$
By (5.2), $\phi_\pi(w{\circ}\xi,w{\circ}\eta)=\phi(\xi,\eta)$ for all $\{\xi,\eta\}\in\Gamma(R^+)$. So $$\phi(\sigma,\tau)\phi(\tau,\mu)+\phi(\sigma,\nu)\phi(\nu,\mu)=0.$$

So $\Gamma(R^+)$ is admissible.
\hfill$\Box$\vspace{3mm}

\section{Diamond (co)chain complex}\vspace{3mm}

\hspace*{5.5mm}{\bf Definition 6.1}
A diamond chain complex $(C,d)$ with basis graph $\Gamma$ is a chain complex satisfying the following conditions.
$\Gamma$ is an admissible diamond graph with sign function $\phi$. $C=\Bbb Z(\Gamma)$ is the abelian group freely generated by $\Gamma$.
For any $a\in\Gamma\subset C$, $da=\Sigma_{|b|=|a|-1,\,\{a,b\}\in\Gamma}\,\phi(a,b)b$.

Dually, regard the dual cochain complex $(C,\delta)$ of $(C,d)$ as the same group with dual differential
$\delta a=\Sigma_{|b|=|a|+1,\,\{a,b\}\in\Gamma}\,\phi(a,b)b$.
Then $(C,\delta)$ is a diamond cochain complex.
\vspace{3mm}

{\bf Remark} There is diamond graph $\Gamma$ that is not gradable, i.e., there is no gradation on $\Gamma$
such that for all $\{a,b\}\in\Gamma$, $|a|=|b|{\pm}1$.
There is also gradable diamond graph $\Gamma$ that is not admissible, i.e., there is no sign function $\phi$
on all edges of $\Gamma$ such that $\phi(a,b)\phi(b,c)+\phi(a,d)\phi(d,c)=0$
for all diamond $\langle a;b;c;d\rangle$.
But if $\Gamma$ is admissible, all sign functions $\phi$ induce isomorphic diamond (co)chain complexes. See \cite{h}.
\vspace{3mm}

{\bf Theorem 6.2}\, {\it $(\Lambda(R^+),d)$ and  $(\Lambda(R^+),\delta)$ are respectively diamond chain and cochain complex with basis graph $\Gamma(R^+)$.
}\vspace{3mm}

{\it Proof}\, Corollary of Theorem 5.2 and Theorem 5.5.
\hfill $\Box$
\vspace{3mm}

{\bf Theorem 6.3} {\it For a connected diamond graph $\Gamma$, all its vertices have the same number of neighbors.
This number is called the rank of the connected diamond graph $\Gamma$. }\vspace{3mm}

{\it Proof}\, If $\Gamma$ has no edge, then it has only one vertex. The rank of $\Gamma$ is $0$.
If $\Gamma$ has only one edge $\{a,b\}$ with two vertices $a$ and $b$, then the rank of $\Gamma$ is $1$.
Suppose $G$ has a vertex $a$ with $n>1$ neighbors.
Let $b$ be a neighbor of $a$ and $A$ and $B$ are respectively the set of neighbors of $a$ and $b$.
For any $v\in A{\setminus}\{b\}$, the three vertices $v,a,b$ determine a unique
vertex $v'$ such that $\langle v;a;b;v'\rangle$ is a diamond.
$v\to v'$ is a 1-1 correspondence from $A{\setminus}\{b\}$ and $B{\setminus}\{a\}$.
So $A$ and $B$ have the same cardinality.
Since $\Gamma$ is connected, all its vertices have the same number of neighbors.
\hfill $\Box$\vspace{3mm}

{\bf Theorem 6.4} {\it Let $(C,d)$ be a diamond chain complex such that the basis graph $\Gamma$ is connected
with rank $r$ and $\Bbb F_p$ be a field of characteristic $p$.
If $r\not\equiv 0\,{\rm mod}\,p$ for $p>0$ and $r\neq 0$ for $p=0$, then\vspace{1mm}\\
\hspace*{40mm}$H_*(C;\Bbb F_p)=0,\,\,H^*(C;\Bbb F_p)=0$.
}\vspace{3mm}

{\it Proof}\,
For $a\in \Gamma$, $\delta(d(a)){+}d(\delta(a))=\Sigma\,\phi(a,b)\phi(b,c)c$,
where the sum is taken over all $b,c\in\Gamma$ such that $\{a,b\}\in \Gamma$ and
$\{b,c\}\in \Gamma$.
Since $\Gamma$ is a diamond graph, we have that for $c\neq a$,
there are unique $b_1,b_2\in \Gamma$ such that $\langle a;b_1;c;b_2\rangle$ is a diamond of $\Gamma$.
Since $\Gamma$ is admissible, we have $\phi(a,b_1)\phi(b_1,c)+\phi(a,b_2)\phi(b_2,c)=0$.
So
$$\delta(d(a)){+}d(\delta(a))=\Sigma_{\{a,b\}\in\Gamma}\,\phi(a,b)\phi(b,a)a=ra.$$
Thus, $d\delta+\delta d=r$. The theorem holds.
\hfill $\Box$\vspace{3mm}

{\bf Definition 6.5} For a weight $\alpha\in\Omega(R^+)$, its rank $r(\alpha)$ is defined to be
the rank of the connected diamond graph $\Gamma(\alpha)$ (by Theorem 5.3).
We have $r(\alpha)=r(w{*}\alpha)$ for all $w\in W$ by Theorem 2.4.\vspace{3mm}

{\bf Theorem 6.6} {\it For a weight $\alpha\in\Omega(R^+)$, $r(\alpha)=0$ if and only if it
satisfies the equivalent conditions in Theorem 2.5.}\vspace{3mm}

{\it Proof}\, $r(\alpha)=0$ if and only if the connected diamond graph $\Gamma(\alpha)$ has only one vertex.
\hfill$\Box$\vspace{3mm}

{\bf Theorem 6.7}\, {\it Let $\Bbb F_p$ be a field of characteristic $p$.
Then
$$H_*(R^+;\Bbb F_0)=\oplus_{w\in W}\,H_*^{w{*}\varrho}(R^+;\Bbb F_0)
=\oplus_{w\in W}\,\Lambda(w{*}\varrho){\otimes}\Bbb F_0,$$
$$H^*(R^+;\Bbb F_0)=\oplus_{w\in W}\,H^*_{w{*}\varrho}(R^+;\Bbb F_0)
=\oplus_{w\in W}\,\Lambda(w{*}\varrho){\otimes}\Bbb F_0.$$
For a prime $p$ and $\alpha\in\Omega(R^+)$, if $r(\alpha)\not\equiv 0\,{\rm mod}\,p$, then\vspace{2mm}\\
\hspace*{38mm}$H_*^\alpha(R^+;\Bbb F_p)=0,\,\,H^*_\alpha(R^+;\Bbb F_p)=0.$
}\vspace{3mm}

{\it Proof}\, Corollary of Theorem 6.6 and Theorem 6.4.
\hfill $\Box$
\vspace{3mm}

{\bf Remark} For $\Bbb F_0=\Bbb C$, Theorem 6.7 is the special case of Kostant Theorem in \cite{k} for the trivial representation.
So Theorem 6.7 is a generalization of Kostant Theorem to the torsion part of the
integral cohomology. \vspace{3mm}

\section{Combinatorial weight}\vspace{3mm}

\hspace*{5.5mm} {\bf Definition 7.1} Let $R^+=A^+_n$ be the positive root system of the simple Lie algebra
over $\Bbb C$ with Dynkin graph $A_n$. Then
$A^+_n=\{e_{i,j}\,|\,0\leqslant i<j\leqslant n\}$
with Lie bracket
$$[e_{i,j},e_{s,t}]=\left\{\begin{array}{cll}
e_{i,t}&&{\rm if}\,\,j=s,\\
-e_{s,j}&&{\rm if}\,\,i=t,\\
0&&{\rm otherwise.}
\end{array}\right.
$$

The combinatorial weight on $\Gamma(A^+_n)$ is defined as follows.
For $\sigma\in\Gamma(A^+_n)$ and $k=0,1,{\cdots},n$, define
$$i_k={\rm number\,of}\,t\,{\rm such\,that}\,e_{k,t}\in\sigma\,+
{\rm number\,of}\,s\,{\rm such\,that}\,e_{s,k}\not\in\sigma.$$
Then $\o\omega(\sigma)=(i_0,i_1,{\cdots},i_n)$ is the combinatorial weight of $\sigma$.
\vspace{3mm}

{\bf Theorem 7.2} {\it The combinatorial weight and weight satisfies
$$\o \omega=\omega+\rho,$$
where $\rho=(\frac n2,{\cdots},\frac n2)\in \Bbb R^{n+1}$.
}\vspace{3mm}

{\it Proof}\, Take $e_{i,j}\in\Bbb R^{n+1}$ such that the $i$-th coordinate is $-1$ and $j$-th
coordinates is $1$ and all other coordinates are $0$. Then
$$\hspace{23mm}\varrho=\frac 12(-n,-(n{-}2),-(n{-}4),{\cdots},(n{-}4),n{-}2,n).\hspace{23mm}$$

So $\o\omega(\emptyset)=\omega(\emptyset)+\rho$. It is easy to check
that for $e_{i,j}\not\in\sigma$,
$$\o\omega(\sigma{\cup}\{e_{i,j}\})=\o\omega(\sigma)-e_{i,j}.$$
Since $\omega$ satisfies the same formula, the theorem holds.
\hfill$\Box$\vspace{3mm}

With Theorem 7.2, we may replace weight by combinatorial weight.
All the main properties of combinatorial weight of $A^+_n$ has been discussed in \cite{h}.
We list some of them (only for cohomology).\vspace{2mm}

(1) There is a combinatorial weight set $\o\Omega(A^+_n)$ such that
$$(\Lambda(A^+_n),\delta)=\oplus_{(i_0,\cdots,i_n)\in\o\Omega(A^+_n)}(\Lambda(i_0,{\cdots},i_n),\delta).$$

(2) For any $(i_0,{\cdots},i_n)\in\o\Omega(A^+_n)$, there is a reflection isomorphism
$$(\Lambda(i_0,i_1,{\cdots},i_{n-1},i_n),\delta)\cong(\Lambda(n{-}i_n,n{-}i_{n-1},{\cdots},n{-}i_1,n{-}i_0),\delta).$$

(3) For any $(i_0,{\cdots},i_n)\in\o\Omega(A^+_n)$, there is a rotation isomorphism
$$(\Lambda(i_0,i_1,{\cdots},i_{n-1},i_n),\delta)\cong(\Lambda(i_n,i_{0},i_1,{\cdots},i_{n-1}),\delta).$$

(4) For any $(i_0,{\cdots},i_n)\in\o\Omega(A^+_n)$, there is a duality isomorphism
$$(\Lambda(i_0,i_1,{\cdots},i_{n-1},i_n),\delta)\cong(\Lambda(n{-}i_0,n{-}i_{1},{\cdots},n{-}i_{n-1},n{-}i_{n}),d).$$

(5) The rank of combinatorial weights satisfies the following formula,
$$r(\cdots,j{+}1,\cdots,i{-}1,\cdots)=r(\cdots,j,\cdots,i,\cdots)+i{-}j{-}1,$$
where the omitted parts of the combinatorial weights are the same.\vspace{2mm}

(6) $r(i_0,{\cdots},i_n)=0$ if and only if it is a permutation of $(0,1,{\cdots},n)$.\vspace{3mm}

The combinatorial weight is defined independent of Weyl group action.
All the proofs in \cite{h} are combinatorial and some of the results can not be generalized.
For example, (2) and (3) above can not be generalized to other positive systems of simple Lie algebras.
All combinatorial weight $(i_0,{\cdots},i_n)\in\o\Omega(A^+_n)$ satisfies $i_0{+}{\cdots}{+}i_n=\frac 12n(n{+}1)$.

Since there is no fraction in coordinates,
the combinatorial weight is much easier to use than weight.
For example,
$H^*_{(i_0,{\cdots},i_{n+1})}(A^+_n;\Bbb F_p)=0$ if the combinatorial weight is a permutation of $(1,1,2,{\cdots},n{-}1,n,n)$
such that $n\not\equiv 0\,{\rm mod}\,p$, since the rank $r(i_0,{\cdots},i_{n+1})=n$.


\vspace{3mm}

\begin{thebibliography}{99}
\bibitem{a} Adams, J.F.,On the structure and application of the Steenrod algebra,
{\it Math.\ Helv.} {\bf 32}.(1958),180--247
\bibitem{v} Dwyer, W.G., Homology of Integral Upper-Triangular Matrices,
{\it Proceedings of the American Mathematical Society.} {\bf
94},(1985)523--528
\bibitem{k} Kostant, B., Lie algebra cohomology and the generalized Borel-Weil theorem,
{\it Ann.\ Math.} {\bf 74.}(1961),329--387
\bibitem{f} Milnor, J., and Moore,J.C.,On the structure of Hopf algebras,
{\it Ann.\ Math.} {\bf 81}.(1965),211--264
\bibitem{j} Maclane,S., Homology,  Springer Verlag, 1963
\bibitem{d} May, J.P., A general algebraic approach to Steenrod algebra,
{\it Lecture Notes in Mathematics.} {\bf 168},153--231
\bibitem{m} May, J.P., The Cohomology of Restricted Lie Algebras and of Hopf Algebras,
{\it Journal of Algebra.} {\bf 3}, (1966)123-145
\bibitem{r}Ravenel, D.C., Complex Cobordism and Stable Homotopy Groups of spheres, Academic Press. INC. (1986)
\bibitem{d} Zheng, Qibing., A New Massey Product on {\bf Ext}
Groups, {\it Journal of Algebra.} {\bf 183}.378--395(1996)
\bibitem{l} Zheng, Qibing., S-module and the New Massey-Product,
{\it Journal of Algebra.} {\bf 190}.487--497 (1997).
\bibitem{h} Zheng, Qibing., Graphs and the (co)homology of Lie algebras,
arxiv;1107.0235v1[math.AT].\end{thebibliography}
\end{document}